\documentclass[a4paper,12pt]{amsart}

\usepackage{amssymb,amscd,amsthm,epsfig,graphicx,amsmath}

\newcommand{\en}{\subset}

 \def\RR{{\mathbb R}}  \def\TT{{\mathbb T}}
   
 \def\ZZ{{\mathbb Z}}

    \def\cU{\mathcal{U}}

\newtheorem*{teo*}{Theorem}
\newtheorem*{prop*}{Proposition}
\newtheorem*{cor*}{Corollary}
\newtheorem*{goal*}{Goal}

\newtheorem*{quest}{Question}

\newtheorem{teo}{Theorem}[section]
\newtheorem{cor}{Corollary}[section]
\newtheorem{lema}{Lemma}[section]
\newtheorem{prop}{Proposition}[section]

\newcommand{\bi}{\begin{itemize}}
\newcommand{\ei}{\end{itemize}}

\theoremstyle{definition} \theoremstyle{remark}
\newtheorem*{defi}{Definition}

\newcommand{\eps}{\varepsilon}
\newcommand{\diam}{\hbox{diam}}

\newcommand{\esbozo}{\vspace{.05in}{\sc\noindent Sketch }}
\newcommand{\dem}{\vspace{.05in}{\sc\noindent Proof.} }

\newcommand{\lqqd}{\par\hfill {$\Box$}}
\newcommand{\finobs}{\par\hfill{$\diamondsuit$} \vspace*{.1in}}

\author[R. Potrie]{Rafael Potrie}
\address{CMAT, Facultad de Ciencias, Universidad de la Rep\'ublica, Uruguay}
\email{rpotrie@cmat.edu.uy}

\title[Jorge Lewowicz and expansive systems]{On the work of Jorge Lewowicz on expansive systems}

\thanks{Expanded version of a talk given by the author in the conference Dynamical Systems in Montevideo held at Montevideo from 13 to 17 of August 2012. The author was partially supported by CSIC group 618/2010.}

\begin{document}

\maketitle

\begin{abstract}
We will try to give an overview of one of the landmark results of
Jorge Lewowicz: his classification of expansive homeomorphisms of
surfaces. The goal will be to present the main ideas with the hope
of giving evidence of the deep and beautiful contributions he made
to dynamical systems. We will avoid being technical and try to
concentrate on the tools introduced by Lewowicz to obtain these
classification results such as Lyapunov functions and the concept
of persistence for dynamical systems. The main contribution that
we will try to focus on is his conceptual framework and approach
to mathematics reflected by the previously mentioned tools and
fundamentally by the delicate interaction between topology and
dynamics of expansive homeomorphisms of surfaces he discovered in
order to establish his result.
\end{abstract}

\begin{flushright}
\emph{The value of a person resides in his major contribution} \\
Arab proverb freely translated\footnote{Translated from a phrase
in the entrance of the Institut du Monde Arab in Paris , France.}.
\end{flushright}

\section{Introduction}

Among the contributions of Lewowicz to mathematics, it is hard to
ignore what I believe to be his major one: The creation of a
school of dynamical systems in Montevideo. This school is also
highly influenced by his way of looking at mathematics which I
hope will be illustrated in this brief note. The main point is
that it is not only the people who work in expansive systems that
has been influenced by him.

As a disclaimer, I mention that I am by far not the best qualified
to write about Lewowicz's work and that this note does not pretend
to be a summary of all of his contributions to mathematics.
However, as a member of the above mentioned school, and having
been strongly influenced by him, I happily accepted this task and
will try to give a panorama of the results of Lewowicz concerning
expansive homeomorphisms. I recommend reading \cite{Martin} for a
global panorama of Lewowicz's contributions.

Let us start with a simple and elementary definition:

\begin{defi}[Expansive homeomorphism] Let $f:M \to M$ be a homeomorphism of a compact metric space $M$.
We say that $f$ is \emph{expansive} if there exists $\alpha>0$
such that given $x \neq y \in M$ there exists $n \in \ZZ$ such
that $d(f^n(x), f^n(y)) \geq \alpha$. The largest possible
constant $\alpha$ is called the \emph{expansivity constant} of $f$
for the metric $d$. \finobs
\end{defi}

It is important to remark that although the definition depends on
the metric, the notion of expansivity is purely topological and
can be stated for general topological spaces by demanding that
points outside the diagonal $\Delta \en M \times M$ escape by
iteration of $f\times f$ from a fixed neighborhood\footnote{In
more technical wording, that $\Delta$ is a locally maximal set for
$f\times f$.} of $\Delta$.

There are many well known examples of expansive homeomorphisms:
for example subshifts of finite type (as well as hyperbolic sets
of Smale's theory) are expansive; and the dynamics restricted to
the minimal set of a Denjoy counter-example is also expansive. We
will focus mainly on other kind of examples, those whose phase
space is a manifold.

Quoting Lewowicz himself in \cite{Lew-BulletinBrazil}

\begin{quote}
(...){\em expansivity means, from the topological point of view,
that any point of the space $M$ has a distinctive dynamical
behavior. Therefore, a stronger interaction between the topology
of $M$ and the dynamics could be expected.}
\end{quote}

Examples of expansive homeomorphisms on manifolds are given by
Anosov and quasi-Anosov diffeomorphisms (see
\cite{Franks-Anosov,FR}) as well as the well known pseudo-Anosov
maps introduced by Thurston (\cite{Thurston}). Of course, products
of expansive homeomorphisms are expansive. In \cite{ObrienReddy}
it is proved that every surface of positive genus admits an
expansive homeomorphism.

We are now ready to state a landmark result of the work of
Lewowicz (\cite{Lew-BulletinBrazil}):

\begin{teo*} There are no expansive homeomorphisms on the two-dimensional sphere $S^2$.
\end{teo*}

This theorem is highly non-trivial, yet, its statement is
completely simple. Let us remark that there is an independent
proof of this result and the rest of the results in
\cite{Lew-BulletinBrazil} by Hiraide (\cite{Hiraide}).

The concrete purpose of this note is to explain the main ideas
behind this result as well as the classification theorem of
expansive homeomorphisms on surfaces obtained by Lewowicz in
\cite{Lew-BulletinBrazil}. Other contributions will be covered by
Ruggiero, specially those concerning geodesic flows and quotient
dynamics (see also \cite{Ruggiero}), of course, both presentations
will have substantial overlap.

It would not be faire to write about Lewowicz's work without
giving motivations for the study of expansive homeomorphisms. We
will review some motivations in the first sections by reviewing
some of Lewowicz's previous work. Other motivations can be found
along the literature, in particular \cite{Lew-Libro} has a chapter
devoted to that.

\section{Lyapunov functions and topological stability}

We start with a quotation from the introduction of
\cite{Lew-Lyapunov} ``{\em This paper contains some results on
topological stability (see [2,3]) that generalize those obtained
in [2] much in the same way as Lyapunov's direct theorem
generalizes the asymptotic stability results of the hyperbolic
case: if at a critical point, the linear part of a vector field
has proper values with negative real parts, the point is
asymptotically stable and the vector field has a quadratic
Lyapunov function; however, asymptotic stability may also be
proved for vector fields with non-hyperbolic linear
approximations, provided they have a Lyapunov function. In a way
this is what we do here, letting Anosov diffeomorphisms play the
role of the hyperbolic critical point and replacing stability by
topological stability; we get this time a class of topological
stable diffeomorphisms wider than the class of Anosov
diffeomorphisms.}''

Lyapunov functions, introduced by Lewowicz in \cite{Lew-Lyapunov}
play the role of a metric, which in the case of expansive
homeomorphisms is a type of adapted metric which allows to
distinguish the stable and unstable parts of the points which are
nearby a given orbit. Other kinds of adapted metrics have been
then proposed (see \cite{Reddy,Fathi}) but we will focus on
Lyapunov functions that are present transversally in much of
Lewowicz's work (and also in some of his students, see
\cite{Groisman-Tesis,Groisman-SinFijos}). One important tool
introduced in order to construct Lyapunov functions is that of
quadratic forms (or infinitesimal Lyapunov functions) which have
had a strong impact in different directions well beyond expansive
systems as we will explain below.

\begin{defi}[Lyapunov function]
A continuous function $V: U \to \RR$ from a neighborhood $U$ of
the diagonal in $M\times M$ is said to be a Lyapunov function for
$f: M \to M$ iff:
\begin{itemize}
\item[-]  $V(x,x)=0$ for every $x\in M$. \item[-] $V(f(x),f(y)) -
V(x,y) > 0$ for every $x\neq y$.
\end{itemize}
\finobs
\end{defi}

It can be seen as a function which ``sees'' the expansivity in one
step. It is proved in \cite{Lew-BulletinBrazil} (Theorem 1.3) that
these functions characterize expansive homeomorphisms (see
\cite{Fathi} for a different approach):

\begin{teo}\label{Teo-ExpSiySoloSiLyap}
A homeomorphism of a compact metric space is expansive if and only
if it admits a Lyapunov function.
\end{teo}

Lyapunov functions also provide a way of establishing topological
stability of diffeomorphisms (see \cite{Walters} for the Anosov
case) which may not be Anosov.

We recall the definition of topological stability. We say that a
homeomorphism $f: M \to M$ of a compact manifold $M$ is
\emph{topologically stable} if there exists $\eps>0$ such that for
every homeomorphism $g: M \to M$ which is at $C^0$-distance
smaller than $\eps$ of $g$ there exists a continuous surjective
map $h: M \to M$ which semiconjugates $f$ and $g$, that is:

$$ f \circ h = h \circ g $$

Thurston's pseudo-Anosov maps (see \cite{ObrienReddy,Thurston}) do
admit Lyapunov functions (see \cite{Lew-Persistence} Lemma 3.4 or
apply the previous theorem\footnote{In fact Lewowicz uses his
construction of a Lyapunov function to obtain an alternative proof
of expansiveness of pseudo-Anosov maps.}), however, they are not
topologically stable: One can make perturbations of a
pseudo-Anosov map making that some points have their orbit going
``across'' the singularities and which will not be shadowed by no
orbit of the pseudo-Anosov map. See \cite{Lew-Persistence} or look
at the figures in \cite{Lew-Libro} (Figure 2 in page 11) or
\cite{Lew-Cermi-Scholarpedia} (Figure 3).

Therefore, the existence of Lyapunov functions alone is not enough
to get topological stability. One must add a new hypothesis which
can be thought of as a weak topological version of hyperbolicity
(see \cite{Lew-Lyapunov} Section 5 for a more general and precise
definition):

\begin{defi}
We say that a Lyapunov function $V: U \to \RR$ is
\emph{non-degenerate} if for every $x\in M$ there exists a
splitting $T_xM =S_x \oplus U_x$ such that if $C_S(x)$ (resp.
$C_U(x)$) is a cone around $S_x$ (resp. $U_x$) then $V(\cdot, x)$
is positive (resp. negative) in $\hat C_S(x)$ (resp. $\hat
C_U(x)$), the projection of $C_S(x)$ (resp. $C_U(x)$) by the
exponential map in a small neighborhood. \finobs
\end{defi}

\begin{figure}[ht]
\begin{center}
\input{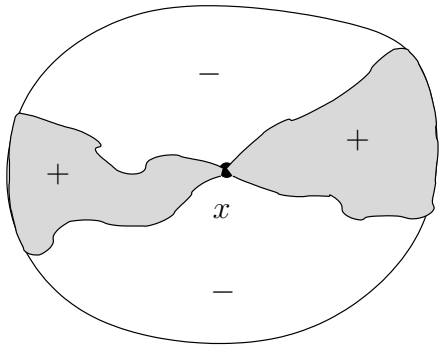}
\caption{\small{Positive and negative regions of $V(\cdot,x)$ when
$V$ is non-degenerate in dimension $2$.}} \label{FiguraLyapunov}
\end{center}
\end{figure}

In a nutshell, the requirement is that the positive and negative
regions of the Lyapunov function in a neighborhood of a point
resemble topologically to the positive and negative part of a
quadratic function (see \footnote{The lines in my drawings are all
crooked on purpose in order to show the topological nature of the
objects, :).} Figure \ref{FiguraLyapunov}). This is exactly what
forbids pseudo-Anosov maps to have non-degenerate Lyapunov
functions in neighborhood of their singularities.

The following theorem is the main theorem of \cite{Lew-Lyapunov}:

\begin{teo}\label{Teo-LyapFunctEntTopSt}
Let $f$ be a $C^1$-diffeomorphism with a non-degenerate Lyapunov
function, then $f$ is topologically stable.
\end{teo}

In \cite{Lew-Lyapunov} a characterization of Anosov
diffeomorphisms in terms of quadratic forms is also given. With
this approach Lewowicz is able to recover classical results on
structural stability of Anosov and characterization of Anosov
systems in terms of cone-families. Quadratic forms turn out to be,
in some applications, better suited for the study of the tangent
map dynamics than cone-fields as we will try to explain in the
next subsection.

\subsection{Quadratic forms, Lyapunov functions and Pesin's theory}

In \cite{Lew-Lyapunov} the following example of diffeomorphism of
$\TT^2$ which is not Anosov and yet admits a non-degenerate
Lyapunov function is proposed:

$$F_c(x,y) = \left(2x - \frac{c}{2\pi} \sin (2\pi x) + y, x-
\frac{c}{2\pi} \sin (2\pi x) +y \right) $$

For $c<1$ the diffeomorphism is Anosov (being linear for $c=0$),
for $c=1$ however, there is no invariant splitting by the
differential in the tangent space of the fixed point $(0,0)$.
However, it can be proved that $F_1$ admits a non-degenerate
Lyapunov function, it is volume preserving and also ergodic.

In \cite{Cats-Enrich} it is proven that $F_1$ as well as many
other examples in the boundary of Anosov diffeomorphisms of
$\TT^2$ are ergodic and non-uniformly hyperbolic. The proof of
non-uniform hyperbolicity relies on the existence of certain
quadratic forms which instead of verifying that their first
difference is everywhere positive, they verify this
almost-everywhere extending the results of \cite{Lew-Lyapunov}.
The following result was stated without a complete proof in
\cite{Lew-Analytic} and the proof was completed in the appendix of
\cite{Markarian-PesinRegion}:

\begin{teo}\label{Teo-FormaCuadratica}
Let $f$ be a volume preserving diffeomorphism admitting a
continuous quadratic form $B: TM \to \RR$ such that the quadratic
form $f^\sharp (B) - B$ is definitely positive almost everywhere.
Then, $f$ is non-uniformly hyperbolic, i.e. Lyapunov exponents are
almost-everywhere non-vanishing.
\end{teo}

Here, we denote $f^\sharp(B)_x(v)= B_{f(x)} (D_x f v)$.

This Theorem was later extended in \cite{Markarian-Toulouse,
Katok} and is quite related to a cone-criterium
(\cite{Wojtkowski}) but works better in some situations (see
\cite{Markarian-PesinRegion,Katok} and references therein). We
will not enter into details about these important results, but we
refer the reader to \cite{Markarian-Libro} for further
developments and applications to billiard systems. Let us just
mention that in the spirit of Lewowicz phrase in the introduction
of \cite{Lew-Lyapunov} and quoted above, the work of
\cite{Markarian-Toulouse} proves a reciprocal statement to the
quadratic form criterium, giving a parallelism with Lyapunov
method and Massera's theorem (an important mentor for Lewowicz),
see \cite{Massera}.

Let us close this section by mentioning a problem which Lewowicz
has always insisted on:

\begin{quest}[Problem 10.3 of \cite{Lew-Cermi-OpenProblems}]
For $c>1$ does the Pesin region of $F_c$ have positive measure?
\end{quest}

The latter is a typical coexistence question which has always
interested many mathematicians. The maps $F_c$ proposed by
Lewowicz are similar to those of \cite{Prytizki} (see also
\cite{Liverani}). See also the work of Pesin (\cite{Pesin}) on the
coexistence problem which is one of the central problems in
dynamics.


\section{Persistence}
\subsection{Persistence vs Topological stability}
The concept of persistency was introduced by Lewowicz in
\cite{Lew-Persistence} in order to study some robust properties of
certain expansive homeomorphisms under perturbations. In a certain
way, it is a property which can be thought of as a dual property
to shadowing.

If an expansive homeomorphism has the shadowing property then it
is topologically stable (see \cite{Lew-Libro}); nevertheless, not
every expansive homeomorphism is topologically stable as we have
already seen. All known expansive homeomorphisms do verify though
this weaker notion of stability which is called \emph{persistence}
(or \emph{semi-persistence}) a term coined by Lewowicz in
\cite{Lew-Persistence} (see also \cite{Lew-Persistence2}).

\begin{defi}[Persistence]
We say that $f: M \to M$ is \emph{persistent} if for every
$\eps>0$ there exists a $C^0$-neighborhood $\cU$ of $f$ such that
for every $g\in \cU$ and $x\in M$ there exists $y\in M$ such that
$$ d(f^n(x),g^n(y)) \leq \eps \qquad \forall n \in \ZZ $$  \finobs
\end{defi}

In Lewowicz words (\cite{Lew-Persistence}):

\begin{quote}``(...)\emph{roughly, the dynamics of $f$ may be found in
each $g$ close to $f$ in the $C^0$-topology; however, these $g$
may present dynamical features with no counterpart in $f$}.''
\end{quote}

In his paper \cite{Lew-Persistence} Lewowicz proves some results
concerning persistence such as persistence for pseudo-Anosov maps
and more generally, for those expansive diffeomorphisms having a
dense set of hyperbolic periodic points with codimension one.
Those results can be thought of as the germ of further
developments in higher dimensions such as
\cite{Vieitez-LocalProdStructure,Vieitez-ClassifWithDensePP,Vieitez-Diffeos,ABP}.

In particular, he shows that a small $C^1$-perturbation of a
pseudo-Anosov map preserving the singularities must be conjugated
to the original map; a kind of structural stability result for
pseudo-Anosov maps. See also \cite{Handel} for further
developments.

Before we continue with the results of \cite{Lew-Persistence} and
some of the consequences found by Lewowicz and coauthors, we are
tempted to add another quote of \cite{Lew-Persistence}:

\begin{quote}``{\em We believe that, apart from such applications, there
is another reason for studying these persistence properties: it
seems plausible to think that if a theory of asymptotic behavior
is possible, then semi-persistence (i.e. persistence of positive
or negative semi-trajectories) should hold on big subsets of $M$
for large classes of dynamical systems}''.
\end{quote}

See \cite{Lew-Persistence2} for advances in that hope. He posed a
precise question about this problem:

\begin{quest}[Problem 10.2 of \cite{Lew-Cermi-OpenProblems}] For an expansive homeomorphism is every semitrajectory
persistent in the future (in the past)?
\end{quest}

Another question which is motivated by this persistency property
can be stated as follows:

\begin{quest} Does an expansive homeomorphism minimize the entropy in its isotopy class?
\end{quest}

This is true for every known example, and it is true after
Lewowicz classification result for surface homeomorphisms (see
also \cite{Handel}). If every expansive homeomorphism is
persistent, then nearby homeomorphisms should have at least the
same topological entropy. However, it seems that the answer of
this fundamental question is at present far out of reach.

\subsection{Expansive systems and stable points}

Probably the first interaction found by Lewowicz between the
topology of the phase space and expansive dynamics is the fact
that a non-trivial compact connected and locally connected set
admitting an expansive homeomorphism cannot have Lyapunov stable
points. If connectedness is not required, this is clearly false as
can be seen by considering an homoclinic orbit between two fixed
points. A more delicate example, where the phase space is
connected can be found in \cite{ReddyRobertson}.

As a way to pave the way of some results in low dimensions which
required hyperbolic periodic points to have codimension one in
dimensions $2$ and $3$, Lewowicz proved in \cite{Lew-Persistence}
the following result:

\begin{teo}[No Stable Points]\label{Teo-NoStable}
Let $f: M \to M$ be an expansive homeomorphism of a non-trivial
compact connected and locally connected metric space, then $f$ has
no Lyapunov stable points.
\end{teo}

Recall that a point $x$ is Lyapunov stable if for every $\eps>0$
there exists $\delta>0$ such that if $d(x,y)< \delta$ then
$d(f^n(x),f^n(y)) < \eps$ for every $n\geq 0$.

We give here a sketch of the proof of this important result:

\esbozo Let $\alpha>0$ be the expansivity constant of $f$. The
proof is divided into 3 steps:

\medskip {\bf Step 1:}  For $\eps< \alpha$, if

$$S_\eps(x) = \{ y \ : \ d(f^n(x),f^n(y))
\leq \eps \ \ n\geq 0 \}$$

\noindent then we have that the diameter of $f^n(S_\eps(x))$
converges to zero uniformly on $x$ and $n$.

\medskip
{\bf Step 2:} If $x$ is a Lyapunov stable point, and $\eps>0$,
there exists $\sigma>0$ such that for $n\geq 0$ we have that
$f^{-n}(S_\eps(x))$ contains the ball of radius $\sigma$ of
$f^{-n}(x)$.

\medskip
{\bf Step 3:} The previous step implies that every point in the
$\alpha$-limit of $x$ is Lyapunov stable. One can prove using this
fact and the first step that the $\alpha$-limit set must consist
of periodic attractors which are only $\alpha$-limit points of
their own orbit. This gives a contradiction, since it implies that
the whole space is a periodic orbit, and being connected a unique
point (contradicting that the space was non-trivial).

\medskip
The hardest step is Step 2 and it is where local connectedness is
used in an essential way. Roughly, using local connectedness, if
the uniform ball cannot be obtained, one finds a sequence of
points $x_n$, $y_n$ such that they are at distance larger than
$\delta$ and remain at distance less than $\eps$ for all future
iterates and for arbitrarily large number of iterates in the past.
Taking limits, one contradicts expansivity.

Let us explain briefly how to find such pair of points: If when
iterating backwards the $\delta$-ball of $x$ there is no uniform
ball, given $n>0$ one can choose an arc $\gamma$ (or a connected
set) with length smaller than $1/n$ and containing $f^{-k_n}(x)$
(where $k_n$ must necessarily tend to $+\infty$ as $n\to +\infty$)
such that $f^{k_n}(\gamma)$ is not contained in $B_\delta(x)$. By
connectedness there exist a point $y_n$ and a backward iterate
$x_n= f^{-m_n}(x)$ at distance larger than $\delta$ and such that
$d(f^j(x_n),f^j(y_n))\leq \eps$ for every $j \geq - k_n + m_n$.

Since the points $x_n$ and $y_n$ are at distance larger than
$\delta$ and its $k_n - m_n$ backward iterate sends them at
distance less than $1/n$ we get that $k_n-m_n$ also goes to
$+\infty$ as $n\to \infty$. Taking convergent subsequences of
$x_n$ and $y_n$ one obtains different points whose orbits remain
at less than $\eps$ for all iterates contradicting expansivity.

\lqqd

\subsection{Analytic models of pseudo-Anosov maps}

In his paper \cite{Lew-Analytic} with E. Lima de Sa, they provide
a new construction of analytic models of pseudo-Anosov maps that
had been obtained by Gerber (\cite{Gerber}) based on previous work
by Gerber with Katok (\cite{GerberKatok}).

The idea is to replace their conditional stability results by the
structural stability theorem of Lewowicz (\cite{Lew-Persistence})
for pseudo-Anosov maps involving the concept of persistence.

It is important to remark that constructing analytic (even smooth)
models of pseudo-Anosov maps is not easy since by a change of
coordinates which is $C^1$ out of a neighborhood of the
singularities one cannot obtain a smooth model (this was shown in
\cite{GerberKatok}), so a more global modification must be made.

The idea involves ``slowing down'' in a neighborhood of the
singularities (much as one does if one wants to smooth the
parametrization of a curve having a corner in its image without
altering the image) and then approximating by analytic maps which
preserve the singularities as well as some $r$-jets of the
derivative of the map in the singularity. This allows to use the
mentioned Lewowicz's results on persistence
(\cite{Lew-Persistence}).

To show how this creation of models is far from being trivial, let
me state an open problem which we are far from understanding. This
question was strongly motivated by discussions with Jorge Lewowicz
and his constant insistence on the lack of understanding we have
of the role of the dynamics of the tangent map (see also the next
subsection for related problems):

\begin{quest}
Let $f: M \to M$ be a \emph{topological Anosov} (i.e. A
homeomorphism of $M$ which preserves two topologically transverse
foliations one of which contracts distances uniformly and the
other one contracts them for backward iterations). Does there
exist a smooth model for $f$? And analytic?. Assuming the previous
questions have positive answers, can these models be made Anosov?.
\finobs
\end{quest}

The question admits a positive answer both in the codimension one
case and in the case where $M$ is a nilmanifold due to the fact
that the classification results of Newhouse-Franks-Manning only
use the fact that the map is a topological Anosov. However, the
question is completely independent a priori of the classification
of Anosov systems.

One of the main contributions of \cite{Lew-Analytic}, though
lateral to the paper has already been explained in this note, and
has to do to the way they prove that the resulting approximation
maps is still Bernoulli with respect to Lebesgue measure (which
can be thought of as the counterpart of the second part of the
question above). To do this, they use quadratic forms and that is
the germ of further results on non-uniform hyperbolicity as we
have already mentioned.

\subsection{The $C^0$-boundary of Anosov diffeomorphisms}

In this section we state a result obtain by Lewowicz in
colaboration with J. Tolosa about the $C^0$-boundary of
codimension one Anosov diffeomorphisms (see \cite{Lew-C0Border}).

They prove:

\begin{teo} Let $f$ be an expansive homeomorphism in the $C^0$-boundary of
Anosov diffeomorphisms of codimension one in $\TT^d$. Then, $f$ is
conjugated to an Anosov.
\end{teo}

With his classification result for expansive homeomorphisms of the
torus, this can be further improved to get:

\begin{teo}
Let $f: \TT^2 \to \TT^2$ be an expansive homeomorphism. Then $f$
is contained in the $C^0$-closure of the set of Anosov
diffeomorphisms of $\TT^2$.
\end{teo}

\dem Consider $h: \TT^2 \to \TT^2$ a homeomorphism isotopic to the
identity such that $f =h \circ A \circ h^{-1}$ where $A$ is a
linear Anosov automorphism. The existence of such an $h$ is given
by Theorem \ref{Teo-Clasificacion} below.

Then there exists a sequence of diffeomorphisms $h_n$ converging
to $h$ in the $C^0$-topology and such that $h_n^{-1}$ also
converges to $h^{-1}$. Since conjugating an Anosov diffeomorphism
by a diffeomorphism gives an Anosov diffeomorphism we get that $f$
is approximated in the $C^0$-topology by Anosov diffeomorphisms.
\lqqd

An important open question that is motivated by this result is the
following:

\begin{quest}[Problem 10.1 of \cite{Lew-Cermi-OpenProblems}]
Does the $C^1$-closure of Anosov diffeomorphisms contains all
expansive diffeomorphisms of $\TT^2$?
\end{quest}

Notice also that Ma\~n\'e has proved that the $C^1$-interior of
expansive diffeomorphisms consists of Quasi-Anosov
ones\footnote{Is in this paper that Ma\~n\'e introduces the
concept of \emph{dominated splitting}.}, in particular in $\TT^2$
of Anosov ones (\cite{Manhe-Expansive}).

Notice that the set of Anosov diffeomorphisms in an given isotopy
class of $\TT^2$ forms a connected set (see \cite{gogolev}).

\section{Classification theorem in surfaces}

It can be shown easily that the only closed one dimensional
manifold, namely the circle, admits no expansive homeomorphisms.
This can be proved using the Poincare's classification of
homeomorphisms of the circle by discussing depending on the
rotation number. Other than that, some examples and some results
on the non-existence of expansive homeomorphisms of other one
dimensional continua, nothing was known about the existence or
structure of expansive homeomorphisms. It is to be remarked that
Ma\~n\'e proved (\cite{Manhe}) that if a compact metric space
admits an expansive homeomorphism, then it must have finite
topological dimension.

Examples in every orientable surface different from the sphere
were already known (\cite{ObrienReddy}), but there was no clue for
example on which isotopy classes admitted them. The classification
of expansive homeomorphisms of surfaces was thus meant to be
started from scratch and that was what Lewowicz did
(\cite{Lew-BulletinBrazil}): He gained an impressive understanding
of their dynamics and their relation with the topology of the
phase space and one of the most striking aspects of his study is
that he relied only on some well known and almost elementary
properties of plane topology. Of course, once he got a
classification of expansive homeomorphisms in terms of their
dynamics and local behavior the final form of the result, giving
conjugation to already known models, used some less elementary
techniques (\cite{Franks-Anosov,Thurston}).

The starting point was the non existence of stable points proved
by him in \cite{Lew-Persistence} and reviewed in the previous
section. In this section we will give an overview of the
classification results for expansive homeomorphisms of surfaces
and the main ideas involved in the proof. We recall that as we
said in the introduction, these results were obtained
independently by Hiraide \cite{Hiraide}.

What we will provide is far from a complete proof of this
classification result, but we hope that the outline here can be
used as a guide to read the original paper
\cite{Lew-BulletinBrazil} and to obtain some insight on the proof.

\subsection{Statement of the result}

Along this section, $S$ will denote an orientable closed (compact,
connected, without boundary) surface. It is well known that these
surfaces are well characterized by their Euler characteristic, and
consist of the sphere $S^2$, the torus $\TT^2$ and the higher
genus surfaces $S_g$ with $g\geq 2$.

The main result of \cite{Lew-BulletinBrazil} is the following:

\begin{teo}[Classification of expansive homeomorphisms of surfaces]\label{Teo-Clasificacion}
Let $f: S \to S$ an expansive homeomorphism. Then, $S \neq S^2$
and:
\begin{itemize}
\item[-] If $S= \TT^2$ then $f$ is conjugate to a linear Anosov
automorphism. \item[-] If $S=S_g$ then $f$ is conjugate to a
pseudo-Anosov map (\cite{Thurston}).
\end{itemize}
\end{teo}

As we mentioned, Lewowicz result has two parts, first, he gives a
complete dynamical classification of expansive homeomorphisms by a
detailed study of the stable and unstable sets of all the points
in $S$, obtaining for them a local product structure outside some
finite set of ``singularities'' which have a local behavior much
like those of pseudo-Anosov maps. Then, by using global arguments
and shadowing results he obtains the desired conjugacy.

Lewowicz result can be though of in a now very fashionable way
called \emph{rigidity}: Rigidity results (or non-existence
results) are those which give strong restrictions from a priori
very weak ones. In the words of Frederic Le Roux in \cite{LeRoux}:
``(...){\em a simple dynamical property can imply a strong
rigidity. The most striking result here is probably
Hiraide-Lewowicz theorem that an expansive homeomorphism on a
compact surface is conjugate to a pseudoAnosov homeomorphism}''.

\subsection{Stable and unstable sets}
This and the next will be the more technical sections of this
note. However, we will try to first give a statement which will be
proved in these two sections which will allow the reader to
continue. Then, we will enter in some details.

Let $f: S \to S$ be an expansive homeomorphism with expansivity
constant equal to $\alpha$. Consider the following sets:

$$ S_\eps(x) = \{ y \in S \ : \ d(f^n(x),f^n(y)) \leq \eps \ \ n\geq 0 \} $$

$$ U_\eps(x) = \{ y \in S \ : \ d(f^{-n}(x),f^{-n}(y)) \leq \eps \ \ n\geq 0 \} $$

Expansivity can be reformulated as $S_\alpha(x) \cap U_\alpha(x) =
\{x\}$ for every $x\in S$. We call $S_\eps(x)$ (resp. $U_\eps(x)$)
the $\eps$-stable set of $x$ (resp. $\eps$-unstable set of $x$).

As we mentioned in the previous section, it can be easily proved
that the diameter of $f^n(S_\eps(x))$ converges to zero uniformly
independently of $x$ if $\eps< \alpha$. The key technical result
in the classification of expansive homeomorphisms of surfaces can
be stated in terms of these sets:

\begin{teo}[Classification Theorem Local Version]\label{Teo-ClasificacionLocal}
Let $f: S \to S$ be an expansive homeomorphism. Then, there exists
a finite set $F$ (possibly empty) such that for every $x\in S
\backslash F$ we have that there exists $\eps>0$ such that
$S_\eps(x)$ is a continuous arc having $x$ in its interior.
Moreover, there exists a neighborhood $U$ of $x$ having local
product structure. For $x \in F$ we have that the sets
$S_\eps(x)\backslash \{x\}$ and $U_\eps(x) \backslash \{x\}$ are
both a finite number ($\geq 3$) of arcs which are alternated and
in each angle they form, there is also local product structure.
\end{teo}

\begin{figure}[ht]
\begin{center}
\input{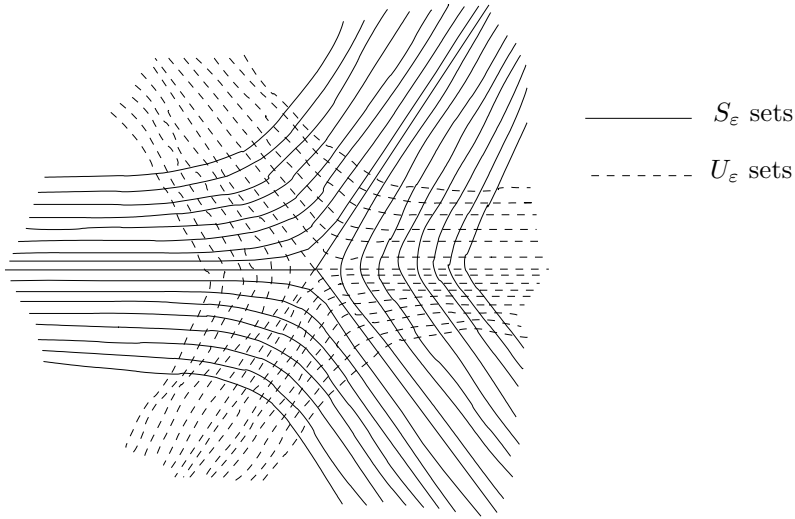} \caption{\small{Local
picture at a singular point with $p=3$ ``legs''.}}
\label{FiguraSingularidad}
\end{center}
\end{figure}

We must explain some of the terminology appearing in the statement
(see also Figure \ref{FiguraSingularidad} for a visual
explanation).

Local product structure means the following: We say that in an
open set $U$ centered in $x$ there is \emph{local product
structure} if there is a homeomorphism
$$h: [-1,1]\times [-1,1] \to \overline U$$ \noindent such that
$h(0,0)=x$ and $h(\{t_0\}\times [-1,1])$ is contained in a stable
set $S_\eps(h(t_0,0))$ and $h([-1,1]\times \{s_0\})$ is contained
in an unstable set $U_\eps(h(0,s_0))$.

In a similar way, given a point $x$, we can consider a connected
component $L^s$ of $S_\eps(x) \backslash \{x\}$ and a connected
component $L^u$ of $U_\eps(x)\backslash \{x\}$. If $U$ is a
neighborhood of $x$ and $A$ is a connected component of $U
\backslash (L^s \cup L^u \cup \{x\})$ which does not intersect
$S_\eps(x) \cup U_\eps(x)$ we say that $A$ is an \emph{angle}. We
say that the angle has local product structure if a similar
property as above holds except that $h: [0,1] \times [0,1] \to
\overline A$ and it sends $h(0,0)=x$ with the rest of the
properties being equal (see Figure \ref{FiguraSingularidad}).

Before we continue with a sketch of the proof of this result, let
us make some comments on some existing extensions. First, similar
properties have been obtained for expansive flows in dimension 3
(\cite{Paternain-FundamentalGroup}). Also, by assuming the
existence of a dense set of topologically hyperbolic periodic
points these results can be extended to any dimension
(\cite{Vieitez-LocalProdStructure,ABP}), except that the behavior
in the singularities is not well understood\footnote{It seems that
we lack examples of ``genuine'' pseudo-Anosov maps in higher
dimensions.} except in dimension 3 or in the codimension one case
(\cite{Vieitez-ClassifWithDensePP,ABP}) where one can show that
they do not exist. With some differentiability assumptions, the
hypothesis of the existence of periodic points can be removed, at
least in dimension 3 (\cite{Vieitez-Diffeos}). This has also been
extended to the plane under certain conditions of the behavior at
infinity (\cite{Groisman-Tesis, Groisman-SinFijos}).

Just to show how far we are from obtaining a similar result in
higher dimensions let me state the following open question (it is
known in the smooth case in $\TT^3$, see \cite{Vieitez-Diffeos}):

\begin{quest} Has every expansive homeomorphism of a manifold a periodic point?
\end{quest}

Now, let us discuss the main points of the proof  of Theorem
\ref{Teo-ClasificacionLocal}. Let us remark that each of this
steps are interesting by themselves, and some of them hold in
higher dimensions.

The first step of the proof consists on showing that every point
in $S$ has a local stable and unstable set of uniform size.

\begin{prop}
For $f: S\to S$ expansive homeomorphism and $\eps<\alpha$ the
expansivity constant, there exists $\delta>0$ such that for every
$x\in S$ we have that the connected component of $S_\eps(x) \cap
B_\delta(x)$ containing $x$ intersects $\partial B_\delta(x)$.
\end{prop}

\esbozo The proof of this proposition holds in any dimension. The
key point is the non existence of Lyapunov stable (in fact,
Lyapunov unstable) points proven in Theorem \ref{Teo-NoStable}.

Once this is obtained, one can construct large connected sets by
considering the sets $D_n$ build as the connected component
containing $x$ of $n$-th preimage of the ball of radius $\eps$
centered at $f^n(x)$ by $f^{-n}$. The fact that there are no
Lyapunov unstable points allows one to prove that these sets have
all diameter bounded from below and allow to construct the desired
set as  $$C_\eps^s(x)= \bigcap_N \overline{\bigcup_{n\geq N}
D_n}$$ One has to check that this has the desired properties (see
\cite{Lew-BulletinBrazil} Lemma 2.1), in particular that the sets
$D_n$ have diameter bounded from bellow. This can be done using
Lyapunov functions and the metric they define, or using the metric
introuduced in \cite{Fathi}. Also, it can be done by barehanded
arguments (see \cite{Lew-Libro}).

\lqqd

We remark that the previous result gives a conceptual proof that
$S^1$ does not admit expansive homeomorphisms: On the one hand
they cannot have Lyapunov stable points, but on the other hand the
stable set of a point must contain a connected set of large
diameter, thus, non-empty interior, a contradiction.

We make a remark on stable and unstable sets which is of
importance in many steps of the proof. It can be thought of as a
``big angles'' result. The proof is not difficult (see Lemma 3.3
of \cite{ABP}).

\begin{prop}[Big Angles]\label{Prop-BigAngles}
Let $f: S\to S$ be an expansive homeomorphism with expansivity
constant $\alpha$. Given $V \en U$ neighborhoods of $x$ and
$\rho>0$ small enough, there exists a neighborhood $W \en V$ of
$x$ such that if $y, z \in W$ we have that $d(S_\eps(y) \cap
U\backslash V, U_\eps(z) \cap U \backslash V) > \rho$.
\end{prop}

The next step of the proof is probably the deepest and it is
really dependent on the two-dimensionality of the problem. Here
one sees a clear manifestation of the already quoted phrase of
``{\em a stronger interaction of the topology of $M$ and the
dynamics of $f$ could be expected} ''.

\begin{teo}\label{Teo-LocalmenteConexo}
For an expansive homeomorphism $f: S \to S$ with expansivity
constant $\alpha$ and $\eps < \alpha/10$, the connected component
of $S_\eps(x)$ containing $x$ is locally connected at each of its
points and therefore arc-connected.
\end{teo}

\esbozo We will only give a brief outline with an heuristic idea
of this subtle proof. We refer the reader to
\cite{Lew-BulletinBrazil} pages 119-121 for details (see also
\cite{Lew-Libro} pages 21-25).

Consider $C_\eps^s(x)$ the connected component of $S_\eps(x)$
containing $x$. We first show that it is locally connected at $x$
and then a clever argument allows to show local connectedness at
every point. Once this is proved, arc-connectedness follows since
a compact connected and locally connected set is arc-connected.

The proof is by contradiction. Roughly, the idea is that if it is
not locally connected at $x$ we can think that in an arbitrarily
small ball of $x$ the set $S_\eps(x)$ is a sequence of connected
sets approaching $x$ but connecting to $C_\eps^s(x)$ outside the
ball. Using separation properties of the plane (which are
extensions of Jordan's curve theorem) we obtain some point $z$
which is trapped in both sides by connected components of
$S_\eps(x)$. Since the unstable set of $z$ has a large connected
component containing $z$, we know it must leave the neighborhood,
however, it can intersect $S_\eps(x)$ only once, so we obtain that
it leaves forming ``small angles'' with $S_\eps(x)$ a
contradiction with Proposition \ref{Prop-BigAngles}.

In fact, there are some subtleties in what we have just said,
since the fact that the unstable set of $z$ has a large connected
component does not imply that it must have two sides, and there is
no problem to have one side going out by intersecting $S_\eps(x)$.
To solve this, Lewowicz makes a clever argument that he then
repeats several times in his proof and so we partially reproduce
it here: He considers an arc joining two different connected
components of $C_\eps^s(x)$ locally and he divides the arc
depending on which side the unstable set of the points leave the
neighborhood: a connectedness argument allows him to conclude that
either there is a point whose unstable intersects twice
$S_\eps(x)$ (contradicting expansivity) or a point whose unstable
leaves forming small angles (also a contradiction). This
connectedness argument uses the fact that stable and unstable sets
vary semicontinuously\footnote{This is a general property that
holds for any homeomorphism and it is not hard to check. See for
example Lemma 3.2 of \cite{ABP}.}.

Now, to get local connectedness at every point, we use local
connectedness at the centers at many scales. Consider $y\in
C_\eps(x) \en C_{2\eps}(y)$. Then $C_{2\eps}(y)$ is locally
connected at $y$, so for every $\sigma>0$ and $z \in C_{\eps}(x)$
close to $y$  there exists a connected set $C \en C_{2\eps}(y)\cap
B_\sigma(y)$ containing $y$ and $z$. Since there are no stable
points, we know that $C \cup C_\eps(x) \en C_{2\eps}(y)$ cannot
separate, so, by an extension of Jordan's separation theorem we
get that $C_\eps(x) \cap C$ is connected and we deduce that
$C_\eps(x)$ is locally connected at $y$.

This finishes the sketch of the proof. \lqqd

We will give an outline of the rest of the proof of Theorem
\ref{Teo-ClasificacionLocal} in the next subsection. We will omit
even more details.

\subsection{Singularities}

The purpose of this section is to outline the rest of the proof of
Theorem \ref{Teo-ClasificacionLocal}. We will not enter in details
here, we will only explain the main steps of the proof.

\begin{figure}[ht]
\begin{center}
\input{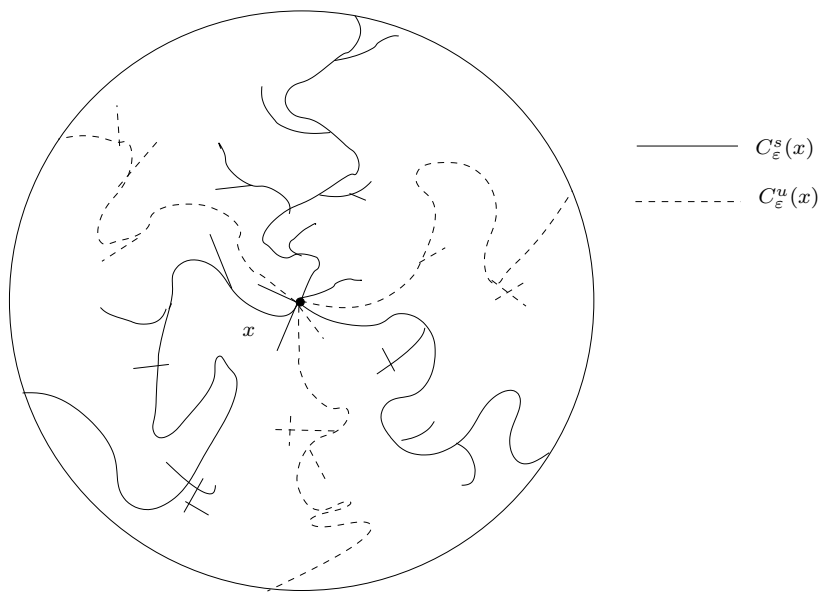} \caption{\small{Structure of
local stable and unstable sets. They are arc connected but may a
priori be still ``ugly''.}} \label{FiguraEstructuraLocal}
\end{center}
\end{figure}

Pick a point $x\in S$. As we have already seen, $S_\eps(x)$ has at
least one connected component intersecting $\partial B_\delta(x)$.
If we consider the connected component $C_\eps^s(x)$ of $S_\eps(x)
\cap B_\delta(x)$ and the connected component $C_\eps^u(x)$ of
$U_\eps(x) \cap B_\delta(x)$ we know that there are arcs joining
$x$ to $\partial B_\delta(x)$ contained in those sets. It is
possible to make an equivalence relation between these arcs that
identify arcs which start at $x$ and then bifurcate near the
boundary of $B_\delta(x)$. By using this, the big angles property
and connected arguments similar to the ones used in the previous
section, Lewowicz shows:

\begin{lema} The number of (equivalence classes of) arcs in $C_\eps^s(x)$ and $C_\eps^u(x)$
joining $x$ to the boundary of $B_\delta(x)$ is the same, finite,
and moreover, they are alternated in the order of $\partial
B_\delta(x)$.
\end{lema}

This result together with the invariance of domain theorem and
further application of the previous arguments give the following
property around points which is almost the end of the proof of
Theorem \ref{Teo-ClasificacionLocal}.

\begin{prop}
For every $x \in S$, there exists a neighborhood $U$ such that
every point $y$ in $U \backslash \{x\}$ has a neighborhood with
local product structure.
\end{prop}

\begin{figure}[ht]
\begin{center}
\input{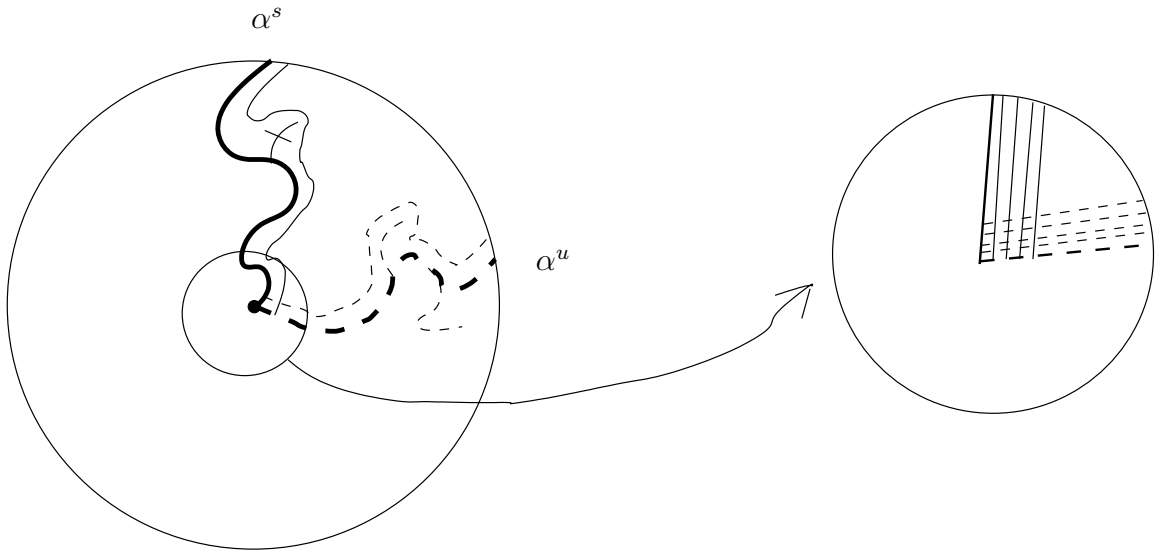} \caption{\small{How to
obtain local product structure.}} \label{FiguraEPL}
\end{center}
\end{figure}

This is proved as follows: Consider an arc $\alpha^s$ of
$C_\eps^s(x)$ and a consecutive one $\alpha^u$ of $C_\eps^u(x)$.
Now, for points of $\alpha^u$ close to $x$ we have that using
semicontinuous variation\footnote{A disclaimer is that to be
precise, this argument needs that there are at least two arcs of
stable and two arcs of unstable for $x$. We will ignore this
problem and ``solve it'' afterwards because we believe it gives a
better heuristic of the global argument. See
\cite{Lew-BulletinBrazil} for a correct proof.} of stable sets and
the ``big angles'' property (Proposition \ref{Prop-BigAngles})
that the stable set of the points near $x$ goes out of
$B_\delta(x)$ near $\alpha^s$. The same happens for points in
$\alpha^s$ near $x$ and their unstable sets. This allows to find a
continuous and injective (due to expansivity) map from a
neighborhood of $x$ in $\alpha^s$ times a neighborhood of $x$ in
$\alpha^u$ into $S$. By the invariance of domain theorem this map
is open and thus every point in this ``angle'' has local product
structure. This can be done in all the angles formed by the stable
and unstable arcs of $x$ (see Figure \ref{FiguraEPL}).

It is immediate to conclude that:

\begin{cor}
There exists a finite set $F \en S$ such that every point outside
of which every point has local product structure. Moreover, for
$x$ in $F$ have a neighborhood such that their local stable and
unstable sets are $p \geq 1$ (and different from $2$ which would
imply local product structure around $x$) arcs starting at $x$ and
arriving at the boundary.
\end{cor}

\begin{figure}[ht]
\begin{center}
\input{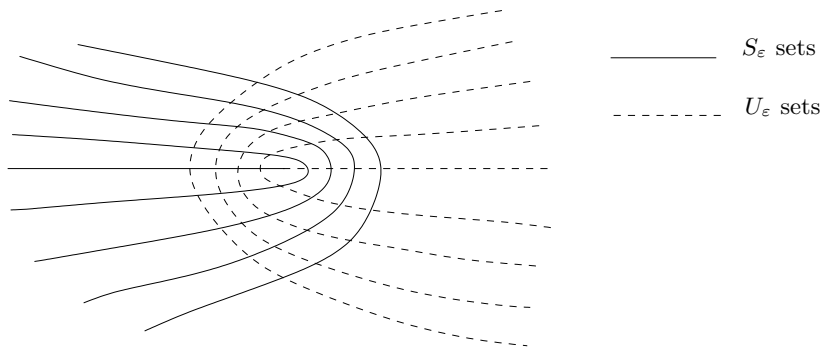} \caption{\small{One leg
implies that local stable and local unstable sets intersect in
more than one point contradicting expansivity.}}
\label{FiguraEspina}
\end{center}
\end{figure}

It remains only to discard the possibility of having a unique arc
in the stable set of $x$. This is an important issue since for
example $S^2$ admits diffeomorphisms (even analytic, see
\cite{Gerber, Lew-Analytic}) that have the local form we have
obtained but with points having a singularity with a unique
``leg''. Needless to say, those examples are not expansive, since
for points very near to $x$ in the stable set, very small
horseshoes are created, contradicting expansivity. Building in
this example, and using the arguments developed by Lewowicz for
the other parts of the proof, one can give a general proof of the
following (see also Figure \ref{FiguraEspina}):

\begin{prop} The number $p$ in the above corollary is $\geq 3$ for every point in $F$.
\end{prop}

This concludes the outline of the proof of Theorem
\ref{Teo-ClasificacionLocal}.

\subsection{Non-existence of expansive homeomorphisms on $S^2$}

We show here how Theorem \ref{Teo-ClasificacionLocal} is enough to
show that the two-dimensional sphere $S^2$ cannot admit expansive
homeomorphisms.

The easiest way to see this is using index theory for foliations.
The local product structure obtained allows one to see that stable
and unstable sets foliate the surface admitting and expansive
homeomorphisms giving rise to a continuous foliation with finitely
many singularities of prong type. Since for every singularity the
number of legs is $\geq 3$ we deduce that even if the foliation
may be non-orientable then the index of the singularities is
always negative (notice that if there were only one leg, then the
index is positive and equal to $1/2$ so that one can make one
example in $S^2$ with four such singularities). This implies that
$S^2$ cannot support such a homeomorphism.

If the reader is not comfortable with the use of continuous (and
not differentiable) foliations, one can go to
\cite{Lew-BulletinBrazil} where a more elementary proof is given
using Poincare-Bendixon's like arguments.

\subsection{Other surfaces}

In the torus case, essentially, due to the work of Franks, it is
enough to show that there are no singularities (which is clear by
the index argument shown above) and that the map is isotopic to a
linear Anosov automorphism. Consider then $f: \TT^2 \to \TT^2$ an
expansive homeomorphism.

Although he might have used the already known argument on the
growth of periodic points and Lefshetz index, Lewowicz gives a
different argument which is very beautiful\footnote{I could not
trace a similar argument to before Lewowicz's paper. However, this
kind of argument has been rediscovered many times so I do not
claim that it is the first time it appeared. I show it in order to
stress the continous search of Lewowicz for understanding and for
conceptual and clean arguments.}.

I outline it here: Lift $f$ to the universal cover to obtain

$$\tilde f: \RR^2 \to \RR^2$$

Let $A \in GL(2,\ZZ)$ be its linear part (i.e. $A$ is the matrix
given by $\tilde f(\cdot) - \tilde f(0) : \ZZ^2 \to \ZZ^2$). If
$A$ is not hyperbolic, then it has both eigenvalues of modulus $1$
since it has determinant of modulus $1$. Then, one obtains that
the diameter of a set iterated by $A$ grows at most polynomially.
Since $\tilde f$ is at bounded distance from $A$, the same holds
for the iterates of a set by $\tilde f$. If $J$ is an unstable arc
contained in a local product structure box, one gets that $\diam
(\tilde f^n(J)) \leq p(n)$ where $p$ is a polynomial.

On the other hand, we know that the length\footnote{Since it is a
continuous arc this is not really well defined. One can measure
thus length by ``counting'' the number of local product structure
boxes it intersects.} of an unstable arc by $\tilde f$ must grow
exponentially due to expansivity, so, for the same $J$ we get that
the length of $\tilde f^n(J)$ is comparable to $\lambda^n$ with
$\lambda>1$. Moreover, since there are no singularities, a
Poincare-Bendixon's like type of argument implies that an arc of
unstable cannot intersect the same box of local product structure
twice. This implies, via the quadratic growth of volume of $\RR^2$
that the diameter of an arc of unstable of length $L$ is
comparable to $\sqrt{L}$ which will still be exponential. This
gives a contradiction and completes the proof. See
\cite{Lew-BulletinBrazil} Theorem 5.3 for more details.

In the higher genus case the proof is even more delicate. He again
stands on previous conjugacy results by Handel \cite{Handel}
(improving the results of \cite{Lew-Persistence}) that state that
in the isotopy class of a pseudo-Anosov map there exist certain
semiconjugacies. Then, as in the torus case he must prove that the
local classification theorem (Theorem
\ref{Teo-ClasificacionLocal}) provides enough tools to show that
$f$ is isotopic to pseudo-Anosov. He uses Thurston's
classification and shows that no homotopy class of simple curves
can be periodic (see Lemma 6.4 of \cite{Lew-BulletinBrazil}) which
allows him to conclude.

\lqqd

%
%
%
\bigskip

{\small \bf Acknowledgments:} {\small My understanding of the
theorem of classification on surfaces owes tremendously to
conversations with A. Artigue and J. Brum. The writing of this
note also benefited from many conversations, exchange and
suggestions of: D.Armentano, E. Catsigeras, M. Cerminara,
M.Delbracio, N. de Leon, H. Enrich, J. Groisman, P. Lessa, R.
Markarian,  A. Passeggi, M. Paternain, R. Ruggiero, A. Sambarino,
M. Sambarino and J. Vieitez. Of course, I want to thank
particularly Jorge Lewowicz, for his work and all the enlightening
conversations we had, mathematically and otherwise. }


\end{document}